\documentclass{sigcomm-alternate}
\usepackage{amsfonts,latexsym,amsmath,amstext,amssymb,verbatim,epsfig,psfrag}
\usepackage{colordvi,pstcol}
\usepackage{graphicx,psboxit}
\usepackage{psfrag}

\def\R{\mbox{I\hspace{-.15em}R}}

\newtheorem{remark}{Remark}
\newtheorem{theorem}{Theorem}

\begin{document}
\title{Note on the equations of diffusion operators associated to a positive matrix}

\numberofauthors{2}
\author{
   \alignauthor Dohy Hong\\
   \affaddr{Alcatel-Lucent Bell Labs}\\
   \affaddr{Route de Villejust}\\
   \affaddr{91620 Nozay, France}\\
   \email{\normalsize dohy.hong@alcatel-lucent.com}
   \alignauthor G\'erard Burnside\\
   \affaddr{Alcatel-Lucent Bell Labs}\\
   \affaddr{Route de Villejust}\\
   \affaddr{91620 Nozay, France}\\
   \email{\normalsize gerard.burnside@alcatel-lucent.com}
}

\date{\today}
\maketitle

\begin{abstract}
In this paper, we describe the general framework to describe the diffusion operators associated to a positive matrix.
We define the equations associated to diffusion operators and present some general properties of their state vectors.
We show how this can be applied to prove and improve the convergence of a fixed point problem associated to the matrix iteration scheme.
The approach can be understood as a decomposition of the matrix-vector product operation in elementary operation
at the vector entry level.
\end{abstract}
\begin{psfrags}
\section{Introduction}\label{sec:intro}
In this paper, we assume that the readers are already familiar with the idea
of the fluid diffusion associated to the D-iteration \cite{d-algo}
to solve the equation:
$$
X = P.X + B
$$
and its application to PageRank equation \cite{dohy}.

For the general description of alternative or existing iteration methods, one
may refer to \cite{Golub1996, Saad}.

\section{Notation}\label{sec:notation}
We will use the following notations:
\begin{itemize}
\item $P \in (\R^+)^{N\times N}$ a positive matrix;
\item $I_d\in (\R^+)^{N\times N}$ the identity matrix;
\item $J_i$ the matrix with all entries equal to zero except for
  the $i$-th diagonal term: $(J_i)_{ii} = 1$;
\item $I = \{i_1,i_2,..,i_n,...\}$ the sequence of nodes for the
  diffusion: $i_k \in \{1,..,N\}$;
\item $\sigma_v : \R^N \to \R$ the scalar product for a given strictly positive vector
   $V>0$: $\sigma_v(X)= <V|X> = \sum_{i=1}^N v_i x_i$;
\item $e$ the normalized unit column vector $1/N (1,..,1)^t$.
\end{itemize}

We say that $P$ is $\sigma_v$-decreasing if:
$$
\forall X\in(\R^+)^N, \sigma_v(PX) \le \sigma_v(X).
$$

We define $P^{\alpha} = (1-\alpha)I_d + \alpha P$.

Then, we have the following results:

\begin{theorem}\label{th:sigma}
$\sigma_v$-decreasing property is stable by composition of operators (matrix product).

If $P$ is $\sigma_v$-decreasing, for $\alpha \ge 0$, $P^{\alpha}$ is $\sigma_v$-decreasing.

If $P$ is $\sigma_v$-decreasing, for $(\alpha,\alpha') \in (\R^+)^2$ such that $\alpha\le\alpha'$,
$\sigma_v(P^{\alpha'}X) \le \sigma_v(P^{\alpha}X)$.
\end{theorem}
\proof
The first point is obvious. The other points are based on the linearity of $\sigma_v$.

\subsection{Diffusion operators}
We define the $N$ diffusion operators associated to $P$ by:
$$
P_i = I_d - J_{i} + P.J_{i}
$$

\begin{theorem}\label{th:Pi}
If $P$ is $\sigma_v$-decreasing, then the diffusion operators
$P_i$ are $\sigma_v$-decreasing. Therefore, for $\alpha \ge 0$, 
$$P^{\alpha}_i = I_d + \alpha(P-I_d).J_i$$ 
is $\sigma_v$-decreasing.
\end{theorem}
\proof
$\sigma_v(P_i.X) = \sigma_v(X) + \sigma_v(P.J_i.X) - \sigma_v(J_i.X)$ and we have
$\sigma_v(P.J_i.X) \le \sigma_v(J_i.X)$, therefore 
$\sigma_v(P_i.X) \le \sigma_v(X)$. The last point is the application of
Theorem \ref{th:sigma} to $P_i$.\\

We recall that the D-iteration is defined by the couple $(P,B) \in \R^{N\times N}\times \R^N$ and exploits two state vectors:
$H_n$ (history) and $F_n$ (residual fluid):
\begin{eqnarray}
F_0 &=& B\\
F_n &=& P_{i_n} F_{n-1}\label{eq:defF}
\end{eqnarray}
and
\begin{eqnarray}\label{eq:defH}
H_0 &=& 0 =(0,..0)^t\\
H_n &=& H_{n-1} + J_{i_n} F_{n-1}.
\end{eqnarray}

The D-iteration is the joint iteration on $(F_n, H_n)$.

Now, we consider a bit more general diffusion iterations as follows:
\begin{eqnarray}
F^{\alpha}_n &=& P^{\alpha_n}_{i_n} F^{\alpha}_{n-1}.\label{eq:defFg}
\end{eqnarray}
and
\begin{eqnarray}\label{eq:defHg}
H^{\alpha}_n &=& H^{\alpha}_{n-1} + \alpha_n J_{i_n} F^{\alpha}_{n-1}
\end{eqnarray}
where $\alpha_n \ge 0$. If for all $n$, $\alpha_n=1$, we have the usual
diffusion iteration.

\begin{theorem}\label{th:HFgen}
$(F^{\alpha}_n, H^{\alpha}_n)$ satisfies:
\begin{align}
H^{\alpha}_n + F^{\alpha}_n &= P. H^{\alpha}_n + B.
\end{align}
\end{theorem}
\proof
The proof is the same as for the $\alpha=1$ (cf. \cite{dohy}) diffusion equations by induction
and using equations \eqref{eq:defFg} and \eqref{eq:defHg}.

\begin{theorem}
Assume we choose $F_0\ge 0$ and $H_0=0$. Then, $F^{\alpha}_n$ and
$H^{\alpha}_n$ are positive and $(H^{\alpha}_n)_i$ is an increasing function for all $i$. 

If $P$ is $\sigma_v$-decreasing, then
$\sigma_v(F^{\alpha}_n)$ is a decreasing function.
\end{theorem}
\proof
The proof is straightforward.

\begin{theorem}\label{th:main}
If we build the two diffusion iterations $(F^{\alpha}_n, H^{\alpha}_n)$ and $(F^{\alpha'}_n,H^{\alpha'}_n)$
from the same initial vector $F_0$ ($H_0=0$) and for the same diffusion sequence $I$,
if for all $n$, $0\le \alpha_n \le \alpha_n'\le 1$, then we have:
\begin{itemize}
\item $\sigma_v(F^{\alpha'}_n) \le \sigma_v(F^{\alpha}_n)$;
\item $H^{\alpha'}_n \ge H^{\alpha}_n$ (for each vector entry);
\item $H^{\alpha'}_n + F^{\alpha'}_n \ge H^{\alpha}_n + F^{\alpha}_n$.
\end{itemize}
\end{theorem}

\proof
The first inequality is a direct consequence of Theorem \ref{th:Pi} and Theorem \ref{th:sigma}.
The third is a direct consequence of the second inequality using Theorem \ref{th:HFgen}.
For the second inequality, we prove by induction: we have obviously $H^{\alpha'}_1 \ge H^{\alpha}_1$.
assume we have, $H^{\alpha'}_n \ge H^{\alpha}_n$. Then, from \eqref{eq:defHg}:
\begin{align*}
H^{\alpha'}_{n+1} &= H^{\alpha'}_{n} + \alpha_{n+1}' J_{i_{n+1}} F^{\alpha'}_{n}\\
&\ge H^{\alpha'}_{n} + \alpha_{n+1} J_{i_{n+1}} F^{\alpha'}_{n}
\end{align*}
and
\begin{align*}
H^{\alpha}_{n+1} &= H^{\alpha}_{n} + \alpha_{n+1} J_{i_{n+1}} F^{\alpha}_{n}
\end{align*}

We need to prove that:
\begin{align*}
H^{\alpha'}_{n} - H^{\alpha}_{n} &\ge \alpha_{n+1} J_{i_{n+1}} (F^{\alpha}_{n}-F^{\alpha'}_n).
\end{align*}
For $i\neq i_{n+1}$, $(H^{\alpha'}_{n} - H^{\alpha}_{n})_i \ge 0$ and $(J_{i_{n+1}} (F^{\alpha}_{n}-F^{\alpha'}_n))_{i} = 0$.
For $i = i_{n+1}$, we only need to handle the case $(F^{\alpha}_{n}-F^{\alpha'}_n)_{i_{n+1}} \ge 0$.
We use the relation: $H^{\alpha'}_{n} - H^{\alpha}_{n} \ge F^{\alpha}_{n}-F^{\alpha'}_n$ to get
the inequality.

\begin{remark}
The power iteration $X_n = P.X_{n-1}$ can be described in the above
scheme  $(F^{\alpha}_n, H^{\alpha}_n)$
with $X_n = F^{\alpha}_{Nn}$
taking $F_0 = X_0$ and
if we apply the cyclic sequence $1,..,N$  ($i_n = n \mbox{ mod } N$) where $\alpha_{kN+i}\le 1$ is chosen
such that we diffuse exactly $(P^k.X_0)_i$
(such a value exists and is less than 1 because after the diffusion of nodes $1,..,i$
the residual fluid on $(i+1)$-th node can only  be increased).
\end{remark}

\subsection{Application to DI$+$}

The D-iteration adaptation to the general case (DI$+$) has been presented in
\cite{dip}.
For the sake of simplicity, we show how the above results apply
to the DI$+$ only for the case when $P$ is an irreducible ergodic positive matrix
having 1 as spectral radius.
First, we take for $v$ the left-eigenvector of $P$ such that $v^t.P = v$.

DI$+$'s idea is to apply the diffusion process to $(P, P.e-e)$.

\begin{theorem}
If we choose the sequence of nodes $I$ such that we only diffuse negative
fluids, then the diffusion applied on $(P, P.e-e)$ converges to $X-e$.
\end{theorem}
\proof
The diffusion from $P.e-e$ can be decomposed as the difference of two
diffusion process $(F^{\alpha'}_n,H^{\alpha'}_n)$ and $(F^{\alpha}_n,H^{\alpha}_n)$ as follows:
we start with $F_0=e$. For the $N$ first diffusions, we choose $i_n=n$ and
\begin{itemize}
\item for $P^{\alpha_n}_{i_n}$, $\alpha_n = 0$;
\item for $P^{\alpha_n'}_{i_n}$, $\alpha_n'$ such that we diffuse exactly
$1/N$ (such a value exists and is less than 1 because after the diffusion of nodes $1,..,i$
the residual fluid on $(i+1)$-th node can only be increased).
\end{itemize}
Then we have:
$F^{\alpha}_N = e$, $H^{\alpha}_N = 0$ and $F^{\alpha'}_N = P.e$, $H^{\alpha'}_N = e$.
Then from the $(N+1)$-th diffusion, we apply exactly the same sequence with $\alpha_n = \alpha_n' = 1$.
From Theorem \ref{th:main}, we have $H_n + e = H^{\alpha'}_n - H^{\alpha}_n \ge 0$ and we have
$\sigma_v(F_n) = \sigma_v(F^{\alpha'}_n - F^{\alpha}_n) = 0$.
If we only diffuse negative fluids, this means that $H_n$ is a decreasing function (per entry).
Since we have $0 \le H_n + e \le e$, $H_n$ is convergent.

It has been shown in \cite{dip} that $|F_n|_v = \sum_i |v_i\times (F_n)_i|$ is a decreasing function.
The convergence of $H_n$ implies of course the convergence to zero of $F_n$.

\begin{remark}
We observed that the implicit strategy  on DI$+$ (diffusion of negative fluids) 
can easily bring a computation time reduction factor above 10.
\end{remark}

\begin{remark}
If we mix the diffusion of positive and negative fluids, there is no guarantee
that the algorithm DI$+$ converges and there is snake-configuration example 
to prove that we may oscillate ($1\to 2$ and $1\to 3$; $2\to 4$; $3\to 5$; $4\to 1$; $5\to 1$).
Now one can easily understand that the optimal sequence may require the diffusion
of positive fluids (in the above snake-configuration, the optimal choice is choosing
once the node $1$ which has a positive fluid assuming a weight of 0.5 on $1\to 2$ and $1\to 3$).
If $P$ is irreducible and ergodic, the author conjectures that we still have the
convergence with $H_n+e\le 1$, if one choose only positive fluids.
The DI$+$ algorithm should converge when $P$ is not ergodic.
The DI$+$ algorithm should even converge when $P$ is not irreducible under the condition
that we only take negative fluids.
\end{remark}


\section{Conclusion}\label{sec:conclusion}
We presented the general equations of the diffusion operators
and the general properties of the associated state vectors with
an illustration of the application to prove the convergence
of DI$+$ with a new sequence $I$ choice strategy.

\end{psfrags}
\bibliographystyle{abbrv}
\bibliography{sigproc}

\end{document}